\documentclass[12pt,reqno]{amsart}

\title[Wiener-Hopf operators]{Multipliers and Wiener-Hopf operators on weighted $L^p$ spaces}

\address{LMAM, Universit\'e de Metz, UMR 7122,Ile du Saulcy 57045, Metz Cedex 1, France.}

\email{petkova@univ-metz.fr}

\usepackage{amsmath} 

\usepackage{amssymb}

\def\squarebox#1{\hbox to #1{\hfill\vbox to #1{\vfill}}}

\newcommand{\F}{{\mathcal F}}

\newcommand{\I}{[-\alpha_1, \:\alpha_0]}

\newcommand{\M}{{\mathcal M}}

\newcommand{\Ms}{{\mathcal W}}

\newcommand{\Z}{{\mathbb Z}}

\newcommand{\RR}{{\mathbb R}}

\newcommand{\Ss}{{\bf S}}

\newcommand{\C}{{\mathbb C}}

\newcommand{\N}{{\mathbb N}}

\newcommand{\w}{{\omega}}

\newcommand{\R}{{\Bbb R}}
\newcommand{\rr}{\mathbb R^+}
\newcommand{\lwa}{L_{\w}^p(\rr)}

\newcommand{\lw}{{L_\omega^p({\Bbb R})}}

\newcommand{\Cc}{{C_c(\R)}}

\newcommand{\cc}{C_c(\R^+)}

\newcommand{\CC}{C_c(\rr)}

\newcommand{\tg}{\theta_\gamma}
\newcommand{\E}{{\bf E}}

\def\bV{\big \Vert}

\renewcommand{\Re}{\mathop{\rm Re}\nolimits}

\renewcommand{\Im}{\mathop{\rm Im}\nolimits}

\def\cat#1{{\mathfrak{#1}}}

\theoremstyle{plain}

\newtheorem{thm}{Theorem}

\newtheorem{lem}{Lemma}

\newtheorem{deff}{Definition}

\theoremstyle{definition}

\vspace{1cm}

\numberwithin{equation}{section}
\author[Violeta Petkova]{Violeta Petkova}
\begin{document}
\maketitle

\begin{abstract}
We study the multipliers $M$ (bounded operators commuting with the translations) on weighted spaces $L_\w^p(\R)$. 
We establish the existence 
of a symbol $\mu_M$ for $M$ and some spectral results for the translations $S_t$ and the multipliers. 
We also study the operators $T$ on the weighted space $L^p_{\omega}(\R^+)$ commuting either with the
 right translations $S_t,\: t \in \R^+$, or left translations $P^+S_{-t}, \: t \in \rr,$ and we establish the existence 
of a symbol $\mu$ of $T$. 
We characterize completely the spectrum $\sigma (S_t)$ of the operator $S_t$ proving that 
$$\sigma (S_t) = \{z \in \C: |z| \leq  e^{t\alpha_0}\},$$
 where
 $\alpha_0$ is the growth bound of $(S_t)_{t\geq 0}$. We obtain a similar result for the spectrum of $(P^+S_{-t}),\: t \geq 0.$ Moreover, for an operator $T$ commuting with $S_t, \: t \geq 0,$ we establish the inclusion 
$\overline{\mu({\mathcal O})}\subset \sigma(T)$, where $\mathcal{O}= \{ z \in \C: \Im z < \alpha_0\}$.

\end{abstract}

\section{Introduction}

Let $E$ be a Banach space of functions on $\R$. For $t\in \R$, define the translation by $t$ on $E$ by
$$S_tf(x)=f(x-t), a.e. ,\:\forall f\in E.$$
We call a multiplier on $E$, every bounded operator on $E$ commuting with $S_t$ for every $t\in \R$.
For the multipliers on a Hilbert space we have the existence of a symbol and some spectral 
results concerning the translations and the multipliers are obtained by using this property of the multipliers (see \cite{V9}, \cite{V8}).
 In the arguments exploited in \cite{V9}, \cite{V8} the spectral mapping theorem of Gearhart \cite{G} for semigroups 
in Hilbert spaces plays an essential role. 

The first purpose of this paper is to extend the main results in \cite{V8}, \cite{V9} 
concerning the existence of the symbol of a multiplier as well as the spectral results in the case where $E$ 
is a weighted $L^p_{\omega}(\R)$ space. 
For general Banach spaces the characterization of the spectrum of the semigroup $V(t) = e^{tG}$ 
by the resolvent of its generator $G$ is much more complicated than for semigroups in Hilbert spaces (see  for instance \cite{LS}). 
In particular, the statements of Lemma 1, 2 and 3 (see Section 2) are
 rather difficult to prove and for general Banach spaces this problem remains open. 
In this paper we restrict our attention to $L^p_{\omega}(\R),\: 1 \leq p < \infty,$ weighted spaces. 
The advantage that we take account is that the semigroup of the translations $(S_t)$ preserves the positive functions.
 For semigroups having this special property in the spaces $L_\w^p(\R)$ 
we have a spectral mapping theorem  (see  \cite{EN}, \cite{We1}, \cite{We2}). 
We obtain Theorems 1-4 for multipliers on $L_\w^p(\R)$ and 
in this work we explain only these parts of the proofs which are based on spectral mapping techniques 
and which are different from the arguments used to establish Theorems 1-4 in the particular case $p=2$ (see for more details \cite{V8}, \cite{V9}). \\

For a Banach space $E$ denote by $E^{'}$ the dual space of $E$. For $f\in E$, $g\in E^{'}$, denote by 
$<f,g>$ the duality. Let $p\geq 1$, and let $\w$ be a weight on $\R$. More precisely, $\w$ is a positive, continuous function such that
$$\sup_{x\in \R}\frac{\w(x+t)}{\w(x)}<+\infty, \forall t\in \R.$$
Let $L^p_\w(\R)$ be the set of measurable functions on $\R$ such that
$$\|f\|_{p, \omega} = \Bigl(\int_\R |f(x)|^p\w(x)^p dx\Bigr)^{1/p}<+\infty,\: 1 \leq p < + \infty.$$

Let $\Cc$ (resp. $\cc$ ) be the space of continuous functions on $\R$ (resp. $\R^+$) with compact support in $\R$ (resp. $\R^+$). Notice that $\Cc$ is dense in $L_\w^p(\R)$.\\
In the following we set $E=L_\w^p(\R)$ and we consider only Banach spaces having this form for $1\leq p<+\infty$. In this case
$$\langle f, g \rangle = \int_{\R}  f(x) \bar{g}(x) \omega^2(x) dx$$
and 
$$|\langle f, g \rangle | \leq  \|f\|_{p, \omega} \|g\|_{q, \omega},\: {\rm for}\: 1 < p < +\infty,$$
where $\frac{1}{p} + \frac{1}{q} = 1.$  For $p = 1$, we have 
$$E' = L^{\infty}_{\omega}(\R) = \{ f\: {\rm is}\:{\rm measurable}:\: |f(x)| \omega(x) < \infty, a. e. \}$$ 
and 
$$\|g\|_{\infty, \omega} = {\rm esssup }\:\{|f(x)|\omega(x),\: x\in \R\}.$$
If $M$ is a multiplier on $E$ then, there exists a distribution $\mu$ such that 
$$Mf=\mu*f,\:\forall f\in \CC.$$
For $\phi\in \CC$, the operator 
$$M_\phi:\lw\ni f\longrightarrow \phi*f$$
is a multiplier on $E$. Introduce 
$$\alpha_0=\lim_{t \to +\infty}\ln \|S_t\|^{\frac{1}{t}},\:
\alpha_1= \lim_{t \to +\infty}\ln \|S_{-t}\|^{\frac{1}{t}}.$$
It is easy to see that $\alpha_1 + \alpha_0 \geq 0.$ Consider 
$$ U= \{z\in \C, \: \Im z\in \I\}.$$
 For an operator $T$ denote by $\rho(T)$ the spectral radius of $T$ and by $\sigma(T)$ the spectrum of $T$.
It is well known that $\rho(S_t)=e^{\alpha_0 t},$ for $t\geq 0$.\\

Given a function $f$ and $a\in \C$, denote by $(f)_a$ the function
$$\R\ni x\longrightarrow f(x)e^{ax}$$
and  denote by $\M$ the algebra of the multipliers on $E$. We note by $\hat{g}$ the Fourier transform of a function $g \in L^2(\R).$
Our first result is a theorem saying that every multiplier on $E$ has a representation by a symbol.
\begin{thm}
Let $M$ be a multiplier on $E$. Then\\
$1$)\:For $a\in \I$, we have $(Mf)_a\in L^2(\RR)$, for every $f \in E$ such that $(f)_a\in L^2(\R).$\\
$2$) For $a \in \I$, there exists a function $\nu_a \in L^{\infty}(\R)$ such that 
$$\widehat{(Mf)_a}(x)=\nu_a(x)\widehat{(f)_a}(x),\:\forall f\in E, \:{\rm with}\:(f)_a\in L^2(\R),\:a.e.$$
Moreover, we have $\|\nu_a\|_{\infty}\leq C \|M\|, \:\forall a \in \I$.\\
$3$)\:\: If $\overset{\circ}{U}\neq \emptyset$, there exists a function $\nu \in {\mathcal H}^{\infty}(\overset{\circ}{U})$ such that
$$\widehat {Mf} (z)=\nu(z){\hat f}(z), \: z \in \overset{\circ}{U},\:\:\forall f \in C_c^\infty(\R),$$
where $\widehat{Mf}(ia+x)=\widehat{(Mf)_a}(x)$, for $a \in ]-\alpha_1, \alpha_0[$, $f \in C_c^\infty(\R)$.
\end{thm} 
The function $\nu$ is called the {\bf symbol} of $M$. 
The above result is similar to that established in \cite{V8}, \cite{V9} and the novelty 
is that we treat Banach spaces $L_\w^p(\R)$ and not only Hilbert spaces. 
Define $\mathcal{A}$ as the closed Banach algebra generated by the operators $M_\phi$, for $\phi\in \Cc$.
Notice that $\mathcal{A}$ is a commutative algebra. Our second result concerns the spectra of $S_t$ and $M \in {\mathcal M}.$
\begin{thm} We have
\begin{equation} \label{eq:1.1}
i)\:\:\sigma (S_t) = \{ z\in \C, \:e^{-\alpha_1 t}\leq |z| \leq e^{\alpha_0 t} \},\:\forall t\in \R.
\end{equation} 
Let $M\in \M$ and let $\mu_M$ be the symbol of $M$.\\
ii) We have
\begin{equation} \label{eq:1.2}
\overline{\mu_M({ U})}\subset \sigma(M).
\end{equation}
iii) If $M\in \mathcal{A} $, then we have
\begin{equation} \label{eq:1.3}
\overline{\mu_M(U)}=\sigma(M).
\end{equation}
\end{thm}

The equality (\ref{eq:1.3}) may be considered as a weak spectral mapping property (see \cite{E}) 
for operators in the Banach algebra ${\mathcal A}.$ On the other hand, it is important to note that if $M\in \M$, but $M \notin \mathcal{A}$, in general we have $\overline{\mu_M({ U})}\neq \sigma(M)$. For the space $E=L^1(\R)$, there exists a counter-example (see section 2 and \cite{E}). Thus the inclusion in (1.2) could be strict.\\

In section 3, we obtain similar results for Wiener-Hopf operators on weighted $L_\w^p(\rr)$ spaces.
In the analysis of  Wiener-Hopf operators some new difficulties appear in comparison with the case of multipliers.

Let $\E$ be a Banach space of functions on $\rr$. 
Let $p\geq 1$ and let $\w$ be a weight on $\R^+$. It means that $\w$ is a positive, continuous function such that
$$0<\inf_{x\geq 0}\frac{\w(x+t)}{\w(x)}\leq\sup_{x\geq 0}\frac{\w(x+t)}{\w(x)}<+\infty, \forall t\in \R^+.$$
Let $L^p_\w(\rr)$ be the set of measurable functions on $\R^+$ such that
$$\int_0^\infty |f(x)|^p\w(x)^p dx<+\infty.$$
 Notice that $\cc$ is dense in $\lwa$. 

Let $P^+$ be the projection from $L^2(\R^{-})\oplus \lwa$ into $\lwa$. 
From now we will denote by $\Ss_a$ the restriction of $S_a$ on $\lwa$ for $a\geq 0$ and, for simplicity, $\Ss_1$ will
 be denoted by $\Ss$. Let $I$ be the identity operator on $\lwa$.
\begin{deff}
A bounded  operator $T$ on $\lwa$ is called a Wiener-Hopf operator if
$$P^+\Ss_{-a}T\Ss_af=Tf,\:\forall a\in \R^+,\:f\in \lwa.$$
\end{deff}
As in \cite{V2} we can show that every Wiener-Hopf operator $T$ has a representation by a convolution.
 More precisely, there exists a distribution $\mu$ such that 
$$Tf=P^+(\mu*f),\:\forall f\in C_c^\infty(\rr).$$
If $\phi\in \Cc$, then the operator 
$$\lwa\ni f\longrightarrow P^+(\phi*f)$$
is a Wiener-Hopf operator and we will denote it by $T_\phi$.
Moreover, we have 
$$(P^+\Ss_{-a}\Ss_a)f=f,\:\forall f\in \lwa,$$
but it is obvious that 
$$(\Ss_aP^+\Ss_{-a})f\neq f,$$
for all $f\in \lwa$ with a support not included in $]a,+\infty[$. 
The fact that $\Ss_a$ is not invertible leads to many difficulties in contrast to the case when we deal with the space 
$\lw.$

Let $\E$ be the space $\lwa$.
As above define
$$\cat{a}_0=\lim_{t \to +\infty}\ln \|\Ss_t\|^{\frac{1}{t}},\: 
\cat{a}_1= \lim_{t \to +\infty}\ln \|\Ss_{-t}\|^{\frac{1}{t}}$$
 and set $J=[-\cat{a}_1, \cat{a}_0]$. The next theorem is similar to Theorem 1.\\

\begin{thm}
Let $a\in J$ and let $T$ be a Wiener-Hopf operator. Then for every $f\in \lwa$ such that $(f)_a\in L^2(\rr)$,
we have \\

\begin{equation} \label{eq:r}
(Tf)_a=P^+\F^{-1}(h_a\widehat{(f)_a})
\end{equation}
with $h_a\in L^\infty(\R)$ and
$$\|h_a\|_\infty\leq C\|T\|,$$
where $C$ is a constant independent of $a$. Moreover, if $\cat{a}_1 + \cat{a}_0 > 0$, the function $h$ defined on 
${\mathcal U}=\{z\in \C:\:\Im z\in J\}$ 
by $h(z)=h_{\Im z}(\Re z)$ is holomorphic on $\overset{\circ}{\mathcal U}$.
\end{thm}
\begin{deff}
The function $h$ defined in Theorem $3$ is called the symbol of $T$. 
\end{deff}

We are able to examine the spectrum of the operators in the space $\Ms$  of bounded operators on $\E$ commuting with $(\Ss_t)_{t\geq 0}$ or $(P^+\Ss_{-t})_{t\geq 0}$.\\

Let $\mathcal{O}=\{z\in \C,\: \Im z< \cat{a}_0\}$ and $\mathcal{V}=\{z\in \C,\: \Im z < \cat{a}_1\}.$

\begin{thm} We have

\begin{equation} \label{eq:1.4}
i)\:\:\sigma (\Ss_t) = \{ z\in \C, \:|z| \leq e^{\cat{a}_0 t} \},\:\forall t>0.
\end{equation} 

\begin{equation} \label{eq:1.5}
ii) \:\:\sigma (P^+\Ss_{-t}) = \{ z\in \C, \:|z| \leq e^{\cat{a}_1 t}\}, \:\forall t>0.
\end{equation}

 Let $T\in \Ms$ and let $\mu_T$ be the symbol of $T$.\\
iii) If $T$ commutes with $\Ss_t,\: \forall t \geq 0,$  then we have
\begin{equation} \label{eq:1.6}
\overline{\mu_T({\mathcal O})}\subset \sigma(T).
\end{equation}
iv) If $T$ commutes with $P^+\Ss_{-t},\: \forall t \geq 0,$  then we have
\begin{equation} \label{eq:1.7}
\overline{\mu_T(\mathcal{V})}\subset \sigma(T).
\end{equation}
\end{thm}

The equalities (\ref{eq:1.4}),(\ref{eq:1.5}) generalize the 
well known results for the spectra of the right and left shifts in the space of sequences $l^2$ (see for instance, \cite{R1}). 
However, our proofs are based heavily on the existence of symbols for Wiener-Hopf operators and having in mind Theorem 3, we follow the arguments in \cite{V10}. \\

In section 4, we obtain a sharp spectral result for Wiener-Hopf operators having the form $T_\phi$ with $\phi \in \Cc$. 
This result is established here for operators in spaces $L_\w^p(\rr)$. It is important to note that even for $p = 2$ and for the Hilbert space $L^2_{\omega}(\R^+)$ our result below is new.

\begin{thm}
 
Let $\phi \in \Cc$. Then\\
i) if $supp\: (\phi)\subset \rr$, we have 
$${\overline{\hat{\phi}({\mathcal O})}}=\sigma(T_\phi).$$
ii) if $supp\: (\phi)\subset \R^-$, we have 
$$\overline{\hat{\phi}(\mathcal{V})}=\sigma(T_\phi).$$

\end{thm}

The above result yields a weak spectral mapping property and can be compared with the equality (1.3) in Theorem 2, however
 the proof is more complicated.

\section{Multipliers on $L^p_\w(\R)$}

Recall that we use the notation $E = L^p_{\omega}(\R)$. We start with the following
\begin{lem} Let  $\lambda \in \C$ be such that $e^\lambda\in \sigma(S)$ and let $\Re \lambda=\alpha_0$. Then
there exists a sequence $(f_{n})_{n\in \N}$ of functions of $E$ and an integer $k \in \Z$  so that  
\begin{equation}\label{eq:2.1}
\lim_{n\to\infty}\bV\Bigl(e^{t A} - e^{(\lambda + 2 \pi k i)t}\Bigr) f_{n}\bV=0, 
\:\forall t \in \R,\:\:\|f_{n}\|=1,\:\forall n\in \N.
\end{equation}
\end{lem}

{\bf Proof.} Let $A$ be the generator of the group $(S_t)_{t\in \R}$.
It is clear that the group $(S_t)_{t\in \R}$ preserves positive functions. Since $E = L^p_{\omega}(\R)$ 
the results of \cite{We1}, \cite{We2} say that the spectral mapping theorem holds and 
$$\sigma(e^{tA})\setminus \{0\} = e^{t\sigma(A)} = \{ e^{t\lambda}:\: \lambda \in \sigma(A)\}.$$ 

In particular, for the spectral bound $s(A)$ of $A$ we get 
$$s(A):= \sup \{ \Re z:\: z \in \sigma(A)\} = \alpha_0.$$
Thus $e^{\lambda} \in \sigma(S) \setminus \{0\} = e^{\sigma(A)}$ yields $\lambda  + 2 \pi k i= \lambda_0 \in \sigma(A)$ for some $k \in \Z.$  On the other hand, $\Re \lambda_0 = \alpha_0$, and we deduce that $\lambda_0$ is on the boundary of the spectrum of $A$. By a well known result, this implies that $\lambda_0$ is in the approximative point spectrum of $A$.

Let $\mu_n$ be a sequence such that $\mu_n\to_{n \to \infty} \lambda_0,\: \Re \mu_n > \lambda_0, \: \forall n \in \N$.
 Then 
$$\bV(\mu_n I - A)^{-1} \bV \geq ({\rm dist}\: (\mu_n , \sigma (A)))^{-1},$$
 hence $\|(\mu_n I - A)^{-1}\| \to  \infty.$
 Applying the uniform boundedness principle and passing to a subsequence of  $\mu_{n}$ (for simplicity also denoted by $\mu_n$),
 we may find $f \in E$ such that 
$$\lim_{n \to \infty} \bV(\mu_{n}I - A)^{-1} f\bV \to \infty.$$
 Introduce $f_{n} \in D(A)$ defined by
$$ f_{n} =  \frac{(\mu_{n}I - A)^{-1} f }{\|(\mu_{n}I - A)^{-1} f \|}.$$
The identity 
$$(\lambda + 2 \pi k i- A) f_{n} = (\lambda_0-\mu_{n}) f_{n} + (\mu_{n} - A) f_{n}$$ 
implies that $(\lambda + 2 \pi k i- A)f_{n} \to 0$ as $n \to \infty.$ 
Then the equality
$$(e^{tA} - e^{t(\lambda + 2 \pi k i)})f_{n} = \Bigl(\int_0^t e^{(\lambda + 2 \pi k i)(t-s)}e^{As} ds\Bigr) (A - \lambda - 2 \pi k i)f_{n}$$
yields (\ref{eq:2.1}). $\Box$

Now we prove the following important lemma.

\begin{lem}
For all $\phi\in C_c^\infty(\R)$ and $\lambda$ such that $e^\lambda \in \sigma (S)$ with $\Re \lambda = \alpha_0$ we have 

\begin{equation}{\label{eq:eqi2}}
|\hat{\phi}(i\lambda +a)|\leq \|M_\phi\|,\:\forall a\in \R.
\end{equation}

\end{lem}

{\bf Proof.}
Let $\lambda\in \C$ be such that $e^\lambda\in \sigma (S)$ and  $\Re \lambda = \alpha_0$ and let $(f_{n})_{n\in \N}$ 
be the sequence constructed in Lemma 1.
We have $$1 = \|f_n\| =\sup_{g\in E^{'},\: \|g\|_{E^{'}} \leq 1}|<f_{n}, g>|.$$
 Then, there exists $g_{n}\in E^{'}$ such that 
$$|<f_{n}, g_{n}>-1|\leq \frac{1}{n}$$
and $\|g_{n}\|_{E^{'}} \leq 1.$
Fix $\phi \in C_c^\infty(\R)$ and consider
$$|\hat{\phi}(i\lambda + a)| \leq  |\hat{\phi}(i \lambda + a) \langle f_n , g_n \rangle| + \frac{1}{n}|\hat{\phi}(i \lambda + a)|$$

$$\leq \Big |\int_{\R} \big \langle \phi(t)\Bigl( e^{(\lambda + 2 \pi k i )t} - S_t\Bigr)
 e^{-i(a + 2 \pi k)t} f_n, g_n \big \rangle dt\Big| + \frac{1}{n}|\hat{\phi}(i \lambda a)|$$
$$+ \Big | \int_{\R} \langle \phi(t) S_t e^{-i(a + 2 \pi k)t} f_n, g_n  \rangle dt\Big |.$$
The first two terms on the right side of the last inequality go to 0 as $n \to \infty$ since by Lemma 1 we have
$$\big \Vert e^{-i(a + 2 \pi k)t}\Bigr(e^{(\lambda + 2 \pi k i)t} - S_t \Bigr)f_n \big \Vert \longrightarrow_{n \to \infty} 0.$$

On the other hand, 
$$I_{n} =\Big |\int_\R <\phi(t) S_t e^{-i(a + 2 \pi k)t} f_{n}, g_{n}> dt \Big|$$
$$ = \Big|\langle \Bigl[ \int_\R \phi(t) e^{-i(a + 2 \pi k)t} f_{n}(. - t)dt\Bigr], g_{n} \rangle\Big|$$
$$= \Big|\langle \int_\R (\phi(. - y) e^{i( a +2 \pi k) y}f_{n}(y) dy, e^{ i( a + 2 \pi k) .}g_{n}\rangle\Big|$$
$$ = \Big|\langle \Bigl(M_{\phi} (e^{i(a + 2 \pi k).} f_{n} )\Bigr), e^{ i(a + 2 \pi k).} g_{n}\rangle\Big|$$
and $I_{n} \leq \|M_{\phi}\|\|f_{n}\| \|g_{n}\|_{E^{'}} \leq  \|M_\phi\|.$ Consequently, we deduce that
$$|\hat{\phi}(i\lambda + a)| \leq \|M_{\phi}\|.$$
$\Box$

Notice that the property (\ref{eq:eqi2}) implies that 
$$|\hat{\phi}(\lambda)|\leq \|M_\phi\|,\:\forall \lambda\in \C,\:{\rm provided} \:\Im \lambda=\alpha_0.$$

\begin{lem}
Let $\phi\in C_c^\infty(\R)$ and let $\lambda$ be such that $e^{-\bar{\lambda}} \in \sigma ((S_{-1})^*)$ with $ \
Re\: \lambda = -\alpha_1$. Then we have
\begin{equation}\label{eq:eqia2}
|\hat{\phi}(i\lambda + a)|\leq \|(M_\phi)\|,\:\forall a\in \R.
\end{equation}

\end{lem}

{\bf Proof.}
Consider the group $(S_{-t})^*_{t\in \R}$ acting on $E^{'}$. Let $\lambda\in \C$ 
be such that $e^{-\bar{\lambda}}\in \sigma (( S_{-1})^*)$ and 
$$|e^{-\bar{\lambda}}|=\rho(S_{-1}) = \rho ((S_{-1})^*) = e^{\alpha_1}.$$
The group $(S_{-t})^*$ preserves
 positive functions. To prove this, assume that $g(x) \geq 0,\: a. e.$ is a positive function and let $h\in E$
be such that $h(x) \geq 0, \: a. e.$ Then
$$\langle h, (S_{-t})^*g \rangle = \langle S_{-t} h, g \rangle \geq 0.$$
If $F(x) =((S_{-t})^*g) (x) < 0$ for $x \in \Lambda \subset \R$ and $\Lambda$ has a positive measure, 
we choose $h(x) = {\bf 1}_{\Lambda}(x).$ Then $\langle {\bf 1}_{\Lambda}, F \rangle \leq 0$ and we conclude that 
$F(x) = 0$ a.e. in $\Lambda$ which is a contradiction.
For the group $(S_{-t})^*$ the spectral mapping theorem holds and,
by the same argument as in Lemma 1, we prove that there exists a sequence $(g_{k})_{k\in \N}$ of functions of $E^{'}$
and an integer $m$ so that for all $t\in \R$,
$$\lim_{k\to\infty}\|(e^{tB}-e^{(-\bar{\lambda} + 2 \pi m i)t)})g_{k}\|_{E'}=0$$
and $\|g_{k}\|_{E^{'}} = 1$.

Since $S_{-t} S_t = I$, we have $(S_t)^*(S_{-t})^* = I$. This implies that 
$$\big \Vert( S_{t})^* g_{k}-e^{(\bar{\lambda} - 2 \pi m i)t} g_{k} \big \Vert_{E'} = \big \Vert\Bigl((S_{t})^* - e^{(\bar{\lambda} - 2 \pi m i)t} (S_t)^*(S_{-t})^*\Bigr) g_{k} \big \Vert_{E'}$$
$$\leq \|(S_t)^*\|_{E'\to E'} \big \Vert \Bigl(e^{(-\bar{\lambda} + 2 \pi m i)t} - (S_{-t})^* \Bigr)g_k \big \Vert_{E'}$$
and we deduce that for every $t \in \R$ we have
$$\lim _{k \to \infty} \big \Vert \Bigl ( (S_t)^* - e^{(\bar{\lambda} - 2 \pi m i)t}\Bigr) g_k\big \Vert_{E'} = 0.$$
For $1 < p < +\infty$ the space $E = L^p_{\omega}(\R)$ is reflexive and the dual to $E'$ can be identified with $E.$ Consequently, since $\|g_k\|_{E'} = 1$, there exists $f_k\in E$ such that 
\begin{equation} \label{eq:2.4}
|<f_k, g_k>-1|\leq \frac{1}{k},\:\: \|f_k\|_{E} \leq 1.
\end{equation}
For $p = 1$ the space $L^1_{\omega}(\R)$ is not reflexive and to arrange (\ref{eq:2.4}), we use another argument. In this case the dual to
$L^1_{\omega}(\R)$ is $L^{\infty}_{\omega}(\R).$ Let $\|g\|_{L^{\infty}_{\omega}(\R)} = 1.$ Fix $0 <\epsilon < 1$ and consider the set
$${\mathcal M}_{\epsilon, m} = \{ x \in \R:\: |g(x)| \omega(x) \geq 1 - \epsilon, \: m \leq x < m + 1\},\: m \in \Z.$$
If $\mu({\mathcal M}_{\epsilon, m})$ (the Lebesgue measure of ${\mathcal M}_{\epsilon, m}$) is zero for all $m \in \Z$, we obtain a contradiction with $\|g\|_{L^{\infty}_{\omega}(\R)} = 1.$ Thus there exists $r \in \Z$ such that $\mu({\mathcal M}_{\epsilon, r}) > 0$.  Now we take
$$f(x) =  \frac{{\bf 1}_{{\mathcal M}_{\epsilon, r}}(x)e^{i \arg (g(x))} }{\mu({\mathcal M}_{\epsilon, r}) \omega^2(x)}.$$
Then 
$$1 \geq \langle f, g\rangle = \int_r^{r + 1}  f(x) \bar{g}(x) \omega^2(x) dx \geq  1 - \epsilon$$
and we can obtain (\ref{eq:2.4}) choosing $\epsilon = 1/k.$
Passing to the proof of (\ref{eq:eqia2}), we get
$$|\hat{\phi}(i\lambda + a) |\leq
 |\hat{\phi}(i\lambda + a)<f_k,g_k>| +\frac{1}{k}|\hat{\phi}(i\lambda + a) |$$

$$ \leq \Big|\int_{\R} < \phi(t) e^{-i(a + 2 \pi m 
)t} f_{k}, \Bigr(e^{\bar{(\lambda}  - 2 \pi m i)t} - ( S_t)^* \Bigr)g_{k} > dt \Big|$$
$$+ \Big|\int_{\R} < \phi(t)S_t\Bigl( e^{-i(a + 2 \pi m)t} f_{k}\Bigr), g_{k} > dt \Big|+\frac{1}{k}|\hat{\phi}(i\lambda + a) |= J_{k}' + I_{k}'+\frac{1}{k}|\hat{\phi}(i\lambda + a) |.$$

 From the argument above we deduce that $J_{k}' \to 0$ as $k \to \infty.$ 
For $I_{k}'$ we apply the same argument as in the proof of Lemma 2 and we deduce 
$$|\hat{\phi}(i\lambda + a)| \leq \|M_{\phi}\|.\:\:\: \Box$$

\vspace{0.4cm}

For the proof of Theorem 1 we apply the argument in \cite{V9} and Lemmas 2-3. 
There exists $ e^{\lambda_0} \in \sigma(S)$ such that $\Re \lambda_0 = \alpha_0.$
 Then for every $z \in \C$ with $\Im z = \alpha_0$ we have
$$|\hat \varphi (z) | \leq \|M_{\varphi}\|.$$
Also there exists $e^{-\lambda_1} \in \sigma((S_{-1})^*)$ with $\Re \lambda_1 = -\alpha_1$
 and for every $z \in \C$ with $\Im z = -\alpha_1$ we have
$$|\hat \varphi (z) | \leq \|M_{\varphi}\|.$$
Applying Phragmen-Lindel\"off theorem for the Fourier transform of $\varphi \in C_c^\infty(\R)$ in the domain $\{z \in \C:\: -\alpha_1 \leq \Im z \leq \alpha_0\}$, we deduce
$$|\hat \varphi (z) | \leq \|M_{\varphi}\|$$
for $z \in U.$ Next we exploit the fact that $M$ can be approximated by $M_{\varphi}$ with respect to the strong operator topology (see
 \cite{V6} for a very general setup covering our case). 
We complete the proof repeating the arguments from \cite{V6}, \cite{V9} and since this leads to minor modifications, we omit the details. To obtain Theorem 2 we follow the same argument as in \cite{V8} and the proof is omitted.

To see that in (\ref{eq:1.2}) the inclusion may be strict, consider  a measure $\eta$ on $\R$ such that 
the operator
$$M_\eta: f\longrightarrow \int_\R S_x(f) d\eta (x) $$
is bounded on $L^1(\R)$. For this it is enough to have $\int_{\R} d|\eta|(x) < \infty.$ 
Then $M_\eta $ is a multiplier on $L^1(\R)$ with symbol  
$$\hat{\eta}(t)=\int_\R e^{-ixt} d\eta(x).$$
On the other hand, there exists a bounded measure $\eta$ on $\R$ such that
$$\overline{\hat{\eta}(\R)}\neq \sigma(M_\eta)$$
 (see for details \cite{E}). 
In $L^1(\R)$ we have $\alpha_0 = \alpha_1 = 0$ and $U=\R$. So we have not the property (\ref{eq:1.2}) in Theorem 2 for every multiplier even in the case $L^1(\R)$.\\

\section{Wiener-Hopf operators}

We need the following lemmas. 

\begin{lem}
 Let $\phi\in C_c(\rr)$. The operator $T_\phi$ commutes with $\Ss_t, \: \forall t > 0,$ if and only 
if the support of $\phi$ is in $\overline{\R^+}$. 
\end{lem}

{\bf Proof.}  Consider $\phi\in C_c(\rr)$ and suppose that $T_{\phi}$ commutes with $\Ss_t,\: t \geq 0.$ We write 
$$\phi = \phi \chi_{\R^{-}} + \phi \chi_{\R^{+}}.$$
If $T_{\phi}$ commutes with $\Ss_t,\: t \geq 0$, then the operator $T_{\phi \chi_{\R^{-}}}$ commutes too. 
Let $\psi = \phi \chi_{\R^{-}}$ and fix $a>0$ such that $\psi$ has a support  in $[-a, 0]$.  
Setting $f=\chi_{[0, a]}$, we get $\Ss_af=\chi_{[a, 2a]}$. 
For $x\geq 0$ we have 
$$P^+(\psi*\Ss_af) (x)= \int_{-a}^0 \psi(t) \chi_{\{a\leq x-t\leq 2a\}}dt
=\int_{max (-a, -2a+x)}^{min(x-a, 0)}\psi(t) dt.$$
Since $P^+(\psi*\Ss_af)=\Ss_a P^+(\psi*f)$, for $x\in [0, a]$, we deduce
$P^+(\psi*\Ss_af) (x)=0$ and 
$$\int_{-a}^{x-a} \psi(t) dt=0,\:\forall x\in [0, a].$$
This implies that 
$\psi(t)=0,$ for $t\in [-a, 0]$ hence $supp(\phi)\subset \overline{\R^+}.$ 

$\Box$

Next we establish the following
\begin{lem}
Let $T_\phi$, $\phi\in C_c(\R)$. Then $T_\phi$ commutes with $P^+(\Ss_{-t})$, $\forall t>0$ if and only if 
$supp(\phi)\subset \overline{\R^-}$. 

\end{lem}

{\bf Proof.}
For $\phi \in C_c(\R)$, suppose that $T_\phi$ commutes with $P^+(\Ss_{-t})$, $\forall t>0$.
 Set $\psi=\phi \chi_{\R^+}$. There exists $a>0$ such that 
$supp(\psi )\subset [0,a].$
We have
$P^+(\psi * P^+ \Ss_{-a} \chi_{[0, a]})=0$ and then $P^+\Ss_{-a}(P^+ \psi *\chi_{[0, a]})=0$. 
This implies that
$$(\psi *\chi_{[0, a]}) (x)=0,\: \forall x>a.$$
On the other hand, for $x > a$  we have
$$(\psi *\chi_{[0, a]})(x)=\int_\R \psi(t) \chi_{[0, a]}(x - t)dt =\int_{max(0, x-a)}^{min(a,x)}\psi(t) dt=\int_{x-a}^a \psi(t) dt.$$
Hence 
$\int_\epsilon^a \psi(t) dt =0, \:\forall a>\epsilon >0$ and we get $\psi=0$. 
Thus we conclude that $supp(\phi)\subset \overline{\R^-}$. $\Box$\\

It is clear that $(\Ss_t)_{t \geq 0}$ and $(P^+ (\Ss_{-t}))_{t \geq 0}$ form continuous semigroups and these semigroups preserve positive functions. Moreover, by using the equality
$$\langle (P^+ \Ss_t)h, g \rangle = \langle h, (P^+ \Ss_{-t})^* g \rangle,$$
we conclude that the semigroup $(P^+ \Ss_{-t})^*$ preserve positive functions.
 The issue is that for $\Ss_t$ and $(P^+ \Ss_{-t})^*$ the spectral mapping theorem holds and we may repeat the  arguments used in section 2. Thus we obtain the following

\begin{lem}
 $1)$ For all $\phi\in C_c^\infty(\R)$ such that $supp(\phi)\subset \R^+$, for $\lambda$ such that $e^\lambda \in \sigma (\Ss)$  and
 $\Re \lambda = {\cat a}_0$, we have 
$$|\hat{\phi}(i\lambda+a)|\leq \|T_\phi\|,\:\forall a\in \R.$$
$2)$ For all $\phi\in C_c^\infty(\R)$ such that $supp(\phi)\subset  \R^-$ and for $\lambda$ such that
 $e^{-\bar{\lambda}} \in \sigma ((\Ss_{-1})^*)$ and $ \
Re\: \lambda = -{\cat a}_1$, we have
$$|\hat{\phi}(i\lambda+a)|\leq \|T_\phi\|,\:\forall a\in \R.$$

\end{lem}


{\bf Proof.}
Let $A$ be the generator of the semi-group $(\Ss_t)_{t\geq 0}$.
First we obtain using the same arguments as in the proof of Lemma 1 that for $\lambda$ such that $e^\lambda\in \sigma(\Ss)$ and 
$\Re \lambda={\cat a}_0$, 
there exists a sequence $(f_n)$ of functions of $\E$ and an integer $k\in \Z$ so that 
$$\lim_{n\to\infty}\bV\Bigl(e^{t A} - e^{(\lambda + 2k \pi  i)t}\Bigr) f_{n}\bV=0, 
\:\forall t \in \R^+,\:\:\|f_{n}\|=1,\:\forall n\in \N.$$
Then we notice that
$$\|(P^+\Ss_{-t}-e^{-(\lambda+2k\pi i)t})f_n\|=\|(P^+\Ss_{-t}-e^{-(\lambda+2k\pi i)t}P^+\Ss_{-t}\Ss_t)f_n\|$$
$$\leq \|P^+\Ss_{-t}\| |e^{-(\lambda+2k\pi i)t}|\|(e^{(\lambda+2k\pi i)t}-\Ss_t)f_n\|, \: \forall t\in \R^+.$$
Thus 
$$\lim_{n\to +\infty}\|(P^+\Ss_{-t}-e^{-(\lambda+2k\pi i)t})f_n\|=0,\:\forall t\in \R^+.$$
So we have
$$\lim_{n\to+\infty}\|\Bigl(P^+\Ss_t - e^{(\lambda + 2 k\pi  i)t}\Bigr) f_{n}\|=0, 
\:\forall t \in \R.$$
 Using the same arguments as in the proof of Lemma 2, we obtain 
$$|\hat{\phi}(i\lambda+a)|\leq \|T_\phi\|,\: \forall a\in \R,\:\forall \phi\in C_c^\infty(\R)$$
and $\lambda$ such that $e^\lambda\in \sigma(\Ss)$ and $\Re \lambda =\cat{a_0}$.
In the same way we prove $2)$ using the semi-group $((P^+\Ss_{-t})^*)_{t\geq 0}$.
$\Box$

To establish  Theorem 3, we use Lemma 6 and we follow with trivial modifications the  arguments in \cite{V2}, \cite{V9}, \cite{V8}. We omit the details. For the proof of Theorem 4 we repeat the arguments in \cite{V10}.\\

Now we pass to the proof of Theorem 5. \\

{\bf Proof of Theorem 5.}

Let $\mathcal{A}$ be the commutative algebra generated by $T_\phi$  for all $\phi$ in $\CC$ with support 
in $\rr$ and $\Ss_x$, for all $x\in \rr$.\\
Denote by $\widehat{\mathcal{A}}$ the set of the characters on $\mathcal{A}$. 
Let $\beta \in \sigma (T_\phi)\setminus \{0\}$. Then there exists $\gamma \in \widehat{\mathcal{A}}$ such that 
$\beta=\gamma(T_\phi).$
We will prove the following equality
$$\gamma(T_\phi)=\int_{\rr} \phi(x)\gamma(S_x) dx.$$
This result il not trivial because we cannot commute $\gamma$ with the Bochner integral $\int_{\rr} \phi(x)S_x dx$. \\

Set
$$\theta_\gamma(x)=\gamma(\Ss_x)=\frac{\gamma(\Ss_x\circ T_\phi)}{\gamma(T_\phi)},\:\forall x\in \rr.$$
Let $\psi\in \CC$ and let $supp (\psi)\subset K$, where $K$ is a compact subset of $\rr$.

Suppose that $(\psi_n)_{n\geq 0}\subset C_K(\rr)$ is a sequence converging to $\psi$ uniformly on $K$.

For every $g\in {\bf E}$, we get
$$\|T_{\psi_n}g-T_\psi g\|\leq\|\psi_n-\psi\|_\infty \sup_{y \in K}\|\Ss_y\| \|g\|$$
and this implies that $\lim_{n \to +\infty}\|T_{\psi_n}-T_\psi\|=0$.  
This shows that the linear map $\psi\longrightarrow T_\psi$
 is sequentially continuous and hence it is continuous from $\CC$ into $\mathcal{A}.$
 Since the map 
$$x\longrightarrow \Ss_x(\psi)$$
is continuous from $\rr$ into $\CC$, we conclude that the map 
$$x\longrightarrow \Ss_x \circ T_\phi =T_{\Ss_x(\phi)}$$ is continuous from $\rr$ into $\mathcal{A}$.
Consequently, the function $\theta_\gamma $ is continuous on $\rr$. 
Introduce 
$$\eta:\CC \ni \psi \longrightarrow \gamma(T_\psi).$$ The map $\eta$ is a continuous linear form on $\CC$ and applying
 Riesz representation theorem, there exists some Borel measure $\mu$ (see for instance, \cite{R}) such that
$$\eta(\psi)=\int_{\rr}\psi(x)d\mu(x), \:\forall \psi\in \CC.$$
This implies that for all $f$, $\psi \in \CC$, we have
$$\gamma(T_\psi \circ T_f)=\int_{\rr} (\psi*f)(t)d\mu(t)$$
$$=\int_{\rr}\Bigl(\int_{\rr} \psi(x)f(t-x)dx\Bigr)d\mu(t).$$
Using the Fubini theorem, we obtain
$$ \gamma(T_\psi \circ T_f)=\int_{\rr} \psi(x) \Big(\int_{\rr} f(t-x)d\mu(t)\Big)dx=
\int_{\rr} \psi(x) \gamma(\Ss_x\circ T_f)dx$$
and replacing $f$  and $\psi$ by $\phi$, we get
\begin{equation}\label{eq:bi}
\gamma(T_\phi)=\int_{\rr}\phi(x)\theta_\gamma(x)dx,\:\forall \phi\in \CC.
\end{equation}

Notice that $\theta_\gamma(x+y)=\theta_\gamma(x)\theta_\gamma(y) ,\:\forall x, y\in \rr$ .
We will prove that $\theta_\gamma(x)\neq 0,\:\forall x\in \rr$.
Suppose $\theta_\gamma(x_0)=0$, for $x_0>0$.
Then $\gamma(\Ss_{x_0})=\Bigl(\gamma(\Ss_{\frac{x_0}{n}})\Bigr)^n=0$ and $\tg(\frac{x_0}{n})=\gamma\Bigl(\Ss_{\frac{x_0}{n}}\Bigr)=0$ for every $n\in \N$. 
Since $\tg$ is continuous on $\rr$, $$\lim_{n \to +\infty} \tg\Bigl(\frac{x_0}{n}\Bigr)=\tg(0)=1$$
and we obtain a contradiction. 
Consequently, we have $\theta_\gamma(x)=\gamma(\Ss_x)\neq 0$, for all $x\in \rr$. 
Now define
$\theta_\gamma(-x)=\frac{1}{\theta_\gamma(x)}, \:\forall x\in \rr.$
It is easy to check that $\tg$ is a  morphism on $\R$.
It is clear that $\tg(x+y)=\tg(x)\tg (y)$, for $(x,y)\in \R^+\times \R^+$ and for $(x,y)\in \R^-\times \R^-$.
Suppose that $x>y>0$, $$\tg(x-y)=\gamma(\Ss_x \Ss_{-y})=\frac{\gamma(\Ss_x \Ss_{-y} \Ss_y)}{\gamma(\Ss_y)}
=\frac{\tg(x)}{\tg(y)}=\tg(x)\tg(-y).$$
Moreover, $$\tg(y-x)=\frac{1}{\tg(x-y)}=\frac{1}{\tg(x)\tg(-y)}=\tg(y)\tg(-x).$$

Since $\tg$ satisfies $\tg(x+y)=\tg(x)\tg(y)$, for all $(x, y)\in \R^2$, it is well known that this implies that
 there exists $\lambda\in\C$ such that 
$\theta_\gamma(x) =e^{\lambda x}$, for all $x\in \R$. 

On the other hand, we have $\gamma(\Ss_x)\in \sigma(\Ss_x)$ and  $\gamma(\Ss_1)=e^{\lambda}\in \sigma(\Ss)$. Thus (\ref{eq:bi}) implies
$$\beta=\gamma(T_\phi)=\hat{\phi}(-i\lambda)$$ with $\lambda \in \mathcal{O}.$
We conclude that 
$$\sigma(T_\phi)\setminus{\{0\}}\subset \hat{\phi}(\mathcal{O}).$$
Now, suppose that $supp(\phi)\subset \R^{-}$. Let $\mathcal{B}$ be the commutative Banach algebra generated by $T_\psi$ for all 
$\psi\in C_c(\R^-)$ and by $P^+\Ss_{-x}$, for all $x\in \rr$.
Let $\kappa\in \sigma(T_\phi)$. Using the same arguments as above, and the set of characters $\widehat{\mathcal{B}}$ of ${\mathcal B}$,
we get 
$$\kappa=\int_{\R^-} \phi(x)e^{\delta x}dx,$$
with $-i\delta\in \mathcal{V}$. This completes the proof of Theorem 5.

$\Box$

\section{Comments and open problems}

Following the general schema of the proof of the existence of symbols for multipliers developed in \cite{V6} for locally compact abelian groups,
it is natural to conjecture that an analog of Theorem 1 holds for general Banach spaces of functions under some hypothesis as we have proved this for general Hilbert space of functions in   \cite{V9}, \cite{V10}.  Using the notations of Section 2, the crucial point is the inequality
\begin{equation} \label{eq:4.1}
|\hat{\varphi}(z)| \leq \|M_{\varphi}\|,\: \forall \varphi \in C_c^{\infty}(\R),\: \Im z = \alpha_0.
\end{equation}
and a similar inequality for $\Im z = - \alpha_1.$
To establish (\ref{eq:4.1}), we introduced the  factor $\langle f_k, g_k \rangle$  (see proof of Lemma 1) close to 1 and we want to estimate 
$\hat{\varphi}(z) \langle f_k, g_k \rangle.$ Here the sequence $f_k, \|f_k\| = 1,$ must be chosen so that for some integers $n_k \in \Z$ and $e^{\lambda} \in \sigma(S),\: \Re \lambda = \alpha_0,$ we have 
\begin{equation} \label{eq:4.2}
\lim_{k \to \infty}\|(S_t - e^{(\lambda + 2 \pi n_k i)t})f_k \| = 0,\:\forall t \in \R.
\end{equation}
If the spectral mapping theorem is true for the group $S_t = e^{At}$, we have $s(A) = \alpha_0$ and (\ref{eq:4.2}) 
can be obtained as in Section 2. On the other hand, if $s(A) < \alpha_0$, we may construct $(f_k)$ assuming that
\begin{equation} \label{eq:4.3}
\sup_{m \in \R} \|(A - \alpha_0 - 2 \pi m i)^{-1}\| = + \infty.
\end{equation}
For Hilbert spaces (\ref{eq:4.3}) holds (see \cite{G}, \cite{LS}, \cite{EN}) and author has
 exploited this property in \cite{V8}, \cite{V9} to complete the proof of (\ref{eq:4.2}). 
For semigroups in Banach spaces $s(A) < \alpha_0$ does not implies in general (\ref{eq:4.3})

 (see a counter-example in Chapter V in  \cite{EN} and the relation between the resolvent of $A$ and 
the spectrum of $S_t$ in \cite{LS}).
 Consequently, it is not possible to use (\ref{eq:4.3}) and to construct a sequence  $f_k$ for which (\ref{eq:4.2}) holds.
 Of course another proof of (\ref{eq:4.1}) could be possible, and in Banach spaces of functions  for which $s(A) < \alpha_0$ this is an open problem.

\end{document}